\documentclass[reqno,a4paper,12pt]{amsart}
\usepackage{amssymb,amsmath,amscd,amstext,amsthm,amsfonts,mathtools}
\usepackage{graphicx}
\usepackage{setspace} 
\usepackage{subfig}
\usepackage[mathscr]{eucal}
\usepackage[ansinew]{inputenc} 

\usepackage{geometry}
\usepackage{amsrefs}
\usepackage{bbm}
\usepackage{float}
\usepackage{caption}
\usepackage{rotating}
\numberwithin{equation}{section}
\newtheorem{theorem}{Theorem}[section]

\newtheorem{corollary}[theorem]{Corollary}
\newtheorem{lemma}[theorem]{Lemma}

{\theoremstyle{definition}
{\newtheorem{remark}[theorem]{Remark}
\newtheorem{example}[theorem]{Example}

\newtheorem{defn}[theorem]{Definition}
}}

\newcommand{\Rr}{{\mathbb{R}}}
\def\diag{\operatorname{diag}}
\newcommand{\matbrackets}[2]{\left[\, {#1} \right]_{{#2}} }
\newcommand{\unit}{\mathbbm{1}}
\renewcommand{\d}{\mathrm d}
\newcommand{\Mat}{{\rm Mat}}

\newtheorem*{theorem*}{Theorem}

\newcommand{\intGamma}{(\Delta^{n-1})^\circ}

\begin{document}
\title[Conservative LV and Replicator Equations]{Conservative Replicator and Lotka-Volterra Equations in the context of Dirac$\backslash$ big-isotropic Structures}
\date{\today}

\author[H. N. Alishah]{Hassan Najafi Alishah}
\address{Departamento de Matem\'atica, Instituto de Ci\^encias Exatas\\
Universidade Federal de Minas Gerais \\
Belo Horizonte, 31270-901, Brazil}
\email{halishah@mat.ufmg.br}

\begin{abstract}
 We introduce an algorithm to find possible constants of motion for a given replicator equation. The algorithm is inspired by Dirac geometry and a Hamiltonian description for the replicator equations with such constants of motion, up to a time re-parametrization, is provided using Dirac$\backslash$big-isotropic structures. Using the equivalence between replicator and Lotka-Volterra (LV) equations, the set of conservative LV equations is enlarged.   Our approach generalizes the well-known use of gauge transformations to skew-symmetrize the interaction matrix of a  LV system. In the case of predator-prey model, our method does allow interaction between different predators and between different preys. 
\end{abstract}

\maketitle
\tableofcontents

\section{Introduction}\label{introduction}
Initially, when Hamiltonian systems appeared in connection with problems of geometric optics and celestial mechanics, the underlying geometry was symplectic geometry. As a compensation for the large number of the properties that Hamiltonian systems have, they are non-generic. Considering other geometries such as poisson , presymplectic, Dirac and big-isotropic helps to extend the scope of the systems that could be described in a Hamiltonian way. For examples of Hamiltonian systems in the context of poisson geometry see references mentioned in \cites{AD2014,MR3223875}, for Hamiltonian systems in the context of presymplectic geometry see references mentioned in \cite{MR3000600} and for Hamiltonian systems described using Dirac structures see \cites{MR3098084,MR2265464}.   

The price to pay for the extension of the scope is losing some of the properties that Hamiltonian systems in the context of symplectic geometry have. But, the main characteristic which is the preservation of Hamiltonian function under the flow of the vector field  is kept. Furthermore, many concepts and techniques, such as completely integrable systems, reduction theory and perturbation theory, which are developed for Hamiltonian systems in the context of symplectic geometry can be carried on for Hamiltonian systems in the context of poisson , presymplectic, Dirac and big-isotropic structures. 

In mathematical biology, replicator equations play a fundamental role in describing evolutionary game dynamics and population dynamics. Consider a population whose individuals interact with each other using one of $n$ possible
pure strategies. 
The state of the population is described by a probability vector $p=(p_1,\ldots, p_n)$, 
with the usage frequency of each pure strategy. This vector is a point in 
 {\em $(n-1)$-dimensional simplex}  
$$\Delta^{n-1}=\{\, (x_1,\ldots, x_n)\in\Rr^ n\,:\, 
x_1+\ldots + x_n=1,\, x_i\geq 0\, \}\;.$$
Given   $x\in\Delta^ {n-1}$, the value
$(A\,x)_i=\sum_{j=1}^n a_{ij}\, x_j$ represents the average pay-off of strategy $i$
within a population at state $x$. Similarly, the value
$x^ t\, A\, x =\sum_{i,j=1}^n a_{ij}\, x_i\,x_j$ stands for the overall average of a population at state $x$,
while the difference $(A\,x)_i - x^ t\, A\, x$ measures the {\em relative fitness} of strategy $i$ in the population $x$.
The {\em replicator} model is the following o.d.e. on $\Delta^{n-1}$
\begin{equation}\label{replicator}
\frac{d x_i}{dt}  = x_i\,\left( (A\,x)_i - x^ t\, A\, x \right) \quad 1\leq i\leq n
\end{equation}
which says that the logarithmic growth rate of each pure strategy's frequency  equals its relative fitness.
The flow of this  o.d.e. is complete and leaves the simplex $\Delta^ {n-1}$ invariant, as well as every of its faces. See \cite{HS} for more details on replicator equations.

Replicator equation~\eqref{replicator} is equivalent to the following equations on $\Rr^{n-1}_+$
\begin{equation}\label{lv-intro}
 \frac{d y_i}{dt} = y_i\left( \, r_i+\sum_{j=1}^{n-1} a^\prime_{ij}\,y_j\right) 
\; (1\leq i\leq n_1))\;,
\end{equation}
where $r_i=a_{in}-a_{nn}$ and $a^\prime_{ij}=a_{ij}-a_{nj}$, see Theorem~\ref{thm:equivalence-lv-replicator}. Equations~\eqref{lv-intro} are known as Lotka-Volterra equations. Lotka-Volterra equations were  introduced independently by Alfred Lotka~\cite{Lot1958} and Vito Volterra \cite{Volt1990} to model the evolution of biological and chemical ecosystems. The affine functions
\[f_i(y)=r_i+\sum_{j=1}^{n-1} a^\prime_{ij}\,y_j\]
are called {\em fitness functions} and $A^\prime=(a^\prime_{ij})\in \Mat_{(n-1)\times (n-1)}(\Rr)$ is called the system's interaction matrix. 

The fitness functions $f_i(y)$ can be more general functions. The dynamics of LV systems can be arbitrarily rich, as was first observed by S. Smale~\cite{Sma1976} who proved that any finite dimensional compact flow can be embedded in a LV system with non linear fitness functions. Later, using a class of embeddings studied by L. Brenig~\cite{Brenig1988}, B. Hern{\'a}ndez-Bermejo and V. Fair{\'e}n ~\cite{BF1997}, it was  proven (see~\cite[Theorems 1 and 2]{BF1997}) that any LV system with polynomial fitness functions can be  embedded in a LV system with affine fitness functions. 
Combining this  with Smale's result, we infer that any finite dimensional compact  flow can be, up to a small perturbation, embedded in a LV system with affine fitness functions. As for applications,  Lotka-Volterra system forms the basis for many models used today in the analysis of population dynamics. It has other applications in Physics, e.g., laser Physics, plasma Physics (as an approximation to the Vlasov-Poisson equation), and neural networks. It appears also in computer science, e.g., communication networks, see \cite{lv-computer}.

In a predator-prey model, the predators will have a negative effect on the prey, and the prey  a positive effect on predator.  Assuming that both interaction effects are equal in size and  ignoring  interactions both between specimens of a single specie and  between different predators or between different  preys gives us a skew-symmetric interaction matrix.  Skew-symmetric interaction matrix yields a conservative Lotka-Volterra equation. In the literature, a LV system with interaction matrix $A^\prime$ is called \emph{conservative} if there exists a diagonal matrix $D^\prime$ with strictly positive element on the diagonal such that $A^\prime D^\prime$ is skew-symmetric. The matrix $A^\prime$ is called \emph{skew-symmetrizable}. Note that a skew-symmetrizable matrix should  necessarily have zero diagonal. The matrix $D^\prime$ was interpreted by Volterra as some sort of normalization by the average weights of the different species.  The change of variable by $D^\prime$ is also called  a gauge transformation for LV system.  If the LV system~\eqref{lv-intro} admits an equilibrium point $q^\prime$, the following function $H:(\Rr^{n-1}_+)^\circ\to\Rr$
\begin{equation}\label{lv-constant-motion-intro}
H(y) = \sum_{j=1}^{n-1} y_j - q^\prime_j\,\log y_j 
\end{equation}
is a constant of motion,  if the system is conservative. Volterra proved that   such conservative LV systems can be given  a Hamiltonian description via symplectic realizatdion and in~\cite{DFO} it is shown that they can be given a Hamiltonian description using a quadratic poisson structure.

{\bf Our results:}
In this paper we consider replicator equations and study their constants of motion and Hamiltonian character using Dirac and big-isotropic structures. Our work here is a continuation of what is done in \cite{AD2014} where for a more general class of equations, i.e.  polymatrix replicator equations, a class of poisson structures on their phase space  was introduced and the corresponding subclass of Hamiltonian polymatrix replicators were characterised. Here (considering only replicator equations and not the general polymatrix replicator) we will be using Dirac and big-isotropic structures instead of poisson structure.

It is shown in \cite{AD2014} that  restriction of a given replicator equation $X_A$ with payoff matrix $A$ and a formal equilibrium point $q=(q_1,\ldots,q_n)$ to the interior of the simplex $\Delta^{n-1}$ is equivalent to a vector field $\tilde{X}_B=B\eta_q$ where $B$ is an $(n-1)\times(n-1)$ matrix and $\eta_q$ is differential of a function defined on $\mathbb{R}^{n-1}$, see Section~\ref{sec:replicator} for more details. Theorem~\ref{thm:consta-of-motion} asserts that if there exist a matrix $D$ such that $DB$ is skew-symmetric and the $1$-form $(1+\sum_{i=1}^{n-1}e^{u_i})\,D^t\,\eta_q$ is closed then $\tilde{X}_B$ and consequently $X_A$ has a constant of motion. The constant of motion for $\tilde{X}_B$ is the primitive of the $1$-form $(1+\sum_{i=1}^{n-1}e^{u_i})\,D^t\,\eta_q$. This result relaxes the skew-symmetric condition on $B$ which is required at \cite{AD2014}. 

In Theorem~\ref{thm:hamiltonian-formalism}, we show that the vector field $Y_B=(1+\sum_{i=1}^{n-1}e^{u_i})\tilde{X}_B,$ is Hamiltonian with respect to a  Dirac$\backslash$big-isotropic which is generated by matrix valued functions ${\bf B}:=(1+\sum_{i=1}^{n-1}e^{u_i}) B$ and ${\bf D}^t:=(1+\sum_{i=1}^{n-1}e^{u_i})D^t$.

The upshot of this approach is detecting new conservative Lotka-Volterra systems via the equivalence between replicator and LV equations.  This result is stated at Theorem~\ref{thm:result-for-lv}. Unlike the case with skew-symmetrizable interaction matrices, conservative LV systems introduced here can have non-zero element on the diagonal. In the case of predator-prey systems, this means the possibility of competition between different predator species or different prey species. Furthermore, instead of requiring an equilibrium, Theorem~\ref{thm:result-for-lv} requires a formal equilibrium which is an equilibrium for the trivial extension of LV equations~\eqref{lv-intro} to $\mathbb{R}^{n-1}$. We also show that the set of conservative LV systems introduced here include the known ones obtained by skew-symmetrization of the interaction matrix $A^\prime$. 

Our results provide a rich set of examples for Hamiltonian systems in the context of symplectic, presymplectic, poisson, Dirac and big-istropic structures. This supports the sentence \emph{"The dynamics of LV systems can be arbitrarily rich"} we stated earlier. 

{\bf Organization of the paper:}
In section~\ref{sec:replicator}, we recall preliminary definitions and statements about replicator equations from~\cite{AD2014} that we need. 
In Section~\ref{algorithm}, we explain the algorithm to find possible constants of motion without any reference to Dirac geometry. In Section~\ref{dirac}, providing a simple introduction to symplectic, poisson , Dirac and big-isotropic structures, we state our results on Hamiltonian formalism of replicator equations. In Section~\ref{lotka-volterra}, the new conservative  Lotka-Volterra equations are been introduced.

\section{Preliminaries}
\label{sec:replicator}
In this section we present preliminary definitions and statements that we will be needing. By 
$$\Delta^{n-1}=\{\, (x_1,\ldots, x_n)\in\Rr^ n\,:\, 
x_1+\ldots + x_n=1,\, x_i\geq 0\, \},\;$$
we denote the {\em $(n-1)$-dimensional simplex}. This simplex is the set of states of a population whose individuals interact with each other using one of the $n$ possible strategies. The {\em replicator} model for the behavior of this population is the following o.d.e. on $\Delta^{n-1}$
\begin{equation}\label{eq:replicator}
\frac{d x_i}{dt}  = x_i\,\left( (A\,x)_i - x^ t\, A\, x \right) \quad 1\leq i\leq n
\end{equation}
where the $n\times n$ matrix $A$ is called the payoff matrix and its entry $a_{ij}$ represents the payoff  
of an individual using pure strategy $i$ when the opponent uses the  pure strategy $j$.
We will denote by $X_A$ the vector filed defined by Equations~\eqref{replicator}. The flow of the vector filed $X_A$ is complete and leaves the simplex $\Delta^{n-1}$ invariant. 
Straightforward calculations show that

\begin{lemma}\label{equivalence} 
\item[1)] The correspondence $A\mapsto X_A$ is linear and its kernel is formed by matrices of the form
\[C=\begin{pmatrix}
c_1&c_2&\ldots&c_n\\c_1&c_2&\ldots&c_n\\\vdots&\vdots&&\vdots\\c_1&c_2&\ldots&c_n
\end{pmatrix}.
\]
Thus, two matrices $A_1,A_2\in\Mat_{n\times n}(\Rr)$ determine the same 
vector field $X_{A_1}=X_{A_2}$ on $\Delta^{n-1}$ iff   the matrix
$A_1-A_2$ has equal rows. 
\item[2)] The projective transformation $\bar{x}_i=\frac{c_ix_i}{\sum_{j=1}^n c_jx_j}$ with $c_j>0$ transforms the vector field $X_A(x)$ to the vector field 
$\frac{1}{(\sum_{j=1}^n c^{-1}_j \bar{x}_j)}. X_{\tilde{A}}(\bar{x})$ where $$\tilde{A}=A\, \begin{pmatrix}
\frac{1}{c_1}&0&\ldots&0\\0&\frac{1}{c_2}&\ldots&0\\\vdots&\vdots&&\vdots\\0&0&\ldots&\frac{1}{c_n}
\end{pmatrix},$$
i.e. to the replicator equation with payoff matrix $\tilde{A}$ up to a a time re-parameterization, see \cite[Exercise $7.1.3$]{HS}. 
\item[3)] Every given face $\sigma$ ($r-$dimensional where $r<n$)  of $\Delta^{n-1}$ is invariant under the replicator equation \eqref{replicator} and the restriction of \eqref{replicator} to $\sigma$ is a replicator equation. 
\end{lemma}

Here we consider replicators equations $X_A$ such that there exists a point $q\in \mathbb{R}^n$ which satisfies $(A\,q)_i=(A \, q)_j$ for all $i,j=1,\ldots,n$. Such a point is an equilibrium point of the natural extension of $X_A$ to the affine subspace generated by $\Delta^{n-1}$. Following \cite{AD2014}, we call such points formal equilibrium. 

Now, we proceed stating some required facts. These facts had been proved  at \cite{AD2014} for  a more general family of equations called polymatrix replicator equations.  We need some notation set up. All the vectors in $\mathbb{R}^n$ will be considered as column vectors and we set $\unit=(1,1,..,1)^t\in\Rr^{n}$ where the dimension will be clear from the context. Also, the letter, $I$, stands for the identity matrix and for every  $x\in \Rr^n$, we denote by $D_x$ the $n\times n$ diagonal matrix $D_x= \diag(x_i)_i$. 

\begin{defn}\label{defn:informal-equilibrium-replicator}
A point $q\in\Delta^{n-1}$ is called a formal equilibrium of $X_{A}$  if  $(A q)_i=(A q)_j$, for $i,j=1,...,n$ and $\sum_{i=1}^nq_i=1$. 
\end{defn}

A formal equilibrium point $q\in\Delta^{n-1}$ is an equilibrium point of $X_A$. Conversely, if $q$ is an  equilibrium point in the interior of the simplex $\Delta^{n-1}$ then it is a formal equilibrium as well, see \cite[Proposition $2$]{AD2014}.

Given a replicator equation $X_A$ let
\begin{align}\label{poissonpmg}
\pi:\mathbb{R}^n&\to \Mat_{n\times n}(\Rr)\\
\pi_A(x):=&(-1)\, T_x\, D_x \, A\, D_x\,T_x^t,\notag
\end{align}
where $T_x= x\, \unit^t -I\;.$ The existence of a formal equilibrium point $q$ yields the fact that $X_A$ can be written in the form
\begin{equation}\label{jacobian-form} X_A(x)=\pi_A(x)\,\d_x H_q  \; \text{ for every }\; x\in\Delta^{n-1}\;,
\end{equation}
where $H_q(x) =\sum_{i=1}^n q_i\log x_i$, see \cite[Proposition $4$]{AD2014}. Notice that the matrix $\pi_A(x)$ defined at  \eqref{poissonpmg} depends on the point $x$. As it is shows in \cite[Theorem $3.5$]{AD2014}, in the interior of the simplex $\Delta^{n-1}$, denoted by $(\Delta^{n-1})^\circ$, one can get ride of this dependence using the diffeomorphism $\phi:\Rr^{n-1}\to \intGamma$ defined by
\begin{equation}\label{phi-alpha}
 \phi(u):= \left(\frac{e^{u_1}}{1+\sum\limits_{i=1}^{{n-1}} e^{u_i}},\ldots ,\frac{e^{u_{{\scriptscriptstyle n-1}}}}{1+\sum\limits_{i=1}^{{n-1}} e^{u_i}},\frac{1}{1+\sum\limits_{i=1}^{{n-1}} e^{u_i}}\right)	\;. 
\end{equation}

For every $u\in\Rr^{n-1}$ and $x=\phi(u)$, we have
\begin{equation}\label{poissonchart}
(\d_u\phi)B(\d_u\phi)^t=(-1)T_xD_xAD_xT^t_x=\pi_A(x)\;,
\end{equation}
where $\d_u \phi$ denotes the jacobian of $\phi$ at point $u$ and  
$B:=(-1)EAE^t $ where $E$ is the  $(n-1)\times n$ matrix

\begin{equation}\label{E-alpha}
E:= \left[\begin{array}{rrrrr}
-1 & 0 & \cdots & 0 & 1 \\
0 & -1 & \cdots & 0 & 1 \\
\vdots & \vdots & \ddots &  \vdots & \vdots \\
0 & 0 & \cdots &  -1 & 1 
\end{array}
\right].
\end{equation}

Denoting the pulled back vector field $(\d_x\phi^{-1})X_A$ by $\tilde{X}_B$ the equality \eqref{jacobian-form} reads as: 
\begin{equation}\label{jacobian-form2}
\tilde{X}_B(u)=B\,\,\d_u(H_q(\phi(u))\quad \forall u\in\mathbb{R}^{n-1}.
\end{equation}

It is clear that any constant of motion for $\tilde{X}_B$ is a constant of motion for $X_A$ in the interior of the simplex $\Delta^{n-1}$. As it is discussed at \cite{AD2014} when $A$, hence $B$, is skew-symmetric then both vector field $X_A$ and $\tilde{X}_B$ are Hamiltonian in the context of poisson geometry, see \cite[Corollary 2.]{AD2014}. A clear consequence is that the function $H(x) =\sum_{i=1}^n q_i\log x_i$ is a constant of motion. 
We will be generalizing this result.


\section{Constants of motion}
\label{algorithm}

In this section, we explain the algorithm to find possible constants of motion without any reference to Dirac geometry.  The starting point is the equality \eqref{jacobian-form2} i.e. 
\[\tilde{X}_B(u)=B\,\,\d_u(H_q(\phi(u))\quad \forall u\in\mathbb{R}^{n-1},\]
where $H_q(\phi(u))=(\sum_{i=1}^{n-1}q_i u_i)-\log(1+\sum_{i=1}^{n-1}e^{u_i})$. We denote by 
\begin{align}\label{eta-formula}\eta_q(u):=d_uH_q(\phi(u))=&\begin{pmatrix}
q_1-\frac{e^{u_1}}{1+\sum_{i=1}^{n-1}e^{u_i}}\\\vdots\\q_{n-1}-\frac{e^{u_{n-1}}}{1+\sum_{i=1}^{n-1}e^{u_i}}
\end{pmatrix}\\
&=\frac{1}{1+\sum_{i=1}^{n-1}e^{u_i}}\begin{pmatrix}
q_1+(q_1-1)e^{u_1}+q_1e^{u_2}+...+q_1e^{u_{n-1}}\\
q_2+q_2e^{u_1}+(q_2-1)e^{u_2}+...+q_2e^{u_{n-1}}\\\vdots\\
q_{n-1}+q_{n-1}e^{u_1}+q_{n-1}e^{u_2}+...+(q_{n-1}-1)e^{u_{n-1}}
\end{pmatrix}\notag
\end{align}

\begin{theorem}\label{thm:consta-of-motion}
Given a replicator equation, $X_A$ with a formal equilibrium point $q$, let $\tilde{X}_B$ be the vector field defined at \eqref{jacobian-form2}. If there exist an $(n-1)\times(n-1)$ matrix $D$ such that 
\begin{enumerate}
\item The matrix $DB$ is anti-symetric.  
\item The $1$-form $(1+\sum_{i=1}^{n-1}e^{u_i})\,D^t\,\eta_q(u)$ is closed, i.e. the matrix 
\[ D^t\begin{pmatrix}
q_1-1&q_1&q_1&\ldots&q_1\\q_2&q_2-1&q_2&\ldots&q_2\\\vdots&\vdots&\vdots&&\vdots\\q_{n-1}&q_{n-1}&q_{n-1}&\ldots&q_{n-1}-1
\end{pmatrix},\]
is diagonal.
\end{enumerate}
 Then the function 
$$H_D(u)=\sum_{i=1}^{n-1}(\sum_{k=1}^{n-1}d_{ki}q_k) u_i+\sum_{i=1}^{n-1}((\sum_{k=1}^{n-1}d_{ki}q_k)-d_{ii})e^{u_i}.$$
 is a constant of motion for the vector field $\tilde{X}_B(u)$ and, consequently, the function 
 $$H_{D}\circ\phi^{-1}(x)=\sum_{i=1}^{n-1}(\sum_{k=1}^{n-1}d_{ki}q_k) \log(\frac{x_i}{x_n})+\sum_{i=1}^{n-1}((\sum_{k=1}^{n-1}d_{ki}q_k)-d_{ii})\frac{x_i}{x_n}$$
  is a constant of motion for the replicator vector field $X_A(x)$.
\end{theorem}
\begin{proof}
By the setting of the theorem $d_uH_D(u)=(1+\sum_{i=1}^{n-1}e^{u_i})\,D^t\,\eta_q(u)$. Also, by equality~\eqref{jacobian-form2} and notation setting~\eqref{eta-formula} we have  
$\tilde{X}_B=B\eta_q(u)$. The following simple observation 
\begin{align*}
2<\tilde{X}_B(u),&d_uH_D(u))>=(1+\sum_{i=1}^{n-1}e^{u_i})\left(<B\eta_q(u),D^t\,\eta_q(u)>+<B\eta_q(u),\,D^t\,\eta_q(u)>\right)\\
&=(1+\sum_{i=1}^{n-1}e^{u_i})<(DB+(DB)^t)\eta_q,\eta_q>(u)=0,
\end{align*}
finishes the proof. 
\end{proof}

\begin{corollary}\label{cor:replicator-general}
Let $q=(q_1,\ldots,q_{n-1},q_n)$ be a point on the affine space generated by $\Delta^{n-1}$ such that $q_n\neq0$,
\begin{equation}\label{Q1}Q_1:=\begin{pmatrix}
q_1-1&q_1&q_1&\ldots&q_1\\q_2&q_2-1&q_2&\ldots&q_2\\\vdots&\vdots&\vdots&&\vdots\\q_{n-1}&q_{n-1}&q_{n-1}&\ldots&q_{n-1}-1
\end{pmatrix},\end{equation}
\begin{equation}\label{Q2}
Q_2:=\frac{1}{q_n}\begin{pmatrix}
q_n+q_1&q_1&\ldots&q_1\\
q_2&q_n+q_2&\ldots&q_2\\
\vdots&\vdots&\vdots\,\,\vdots\,\,\vdots&\vdots\\
q_{n-1}&q_{n-1}&\ldots&q_n+q_{n-1}
\end{pmatrix},
\end{equation}
and $X_A$ be a replicator equation where its payoff matrix $A$ is of the form
\[A=\begin{pmatrix}
a_{11}&a_{12}&\ldots&a_{1(n-1)}&-\sum_{j=1}^{n-1}\frac{q_j}{q_n}a_{1j}\\
a_{21}&a_{22}&\ldots&a_{2(n-1)}&-\sum_{j=1}^{n-1}\frac{q_j}{q_n}a_{2j}\\
\vdots&\vdots&\vdots\,\,\vdots\,\,\vdots&\vdots&\vdots\\
a_{(n-1)1}&a_{(n-1)2}&\ldots&a_{(n-1)(n-1)}&-\sum_{j=1}^{n-1}\frac{q_j}{q_n}a_{(n-1)j}\\
0&0&\ldots&0&0
\end{pmatrix}.\]
Furthermore, let $\bar{Q_1}$ be a matrix such that $(\bar{Q_1})^tQ_1$ is diagonal. If there exists  $d=(d_1,d_2,\dots,d_{n-1})\in\mathbb{R}^{n-1}$ such that the matrix 
\begin{equation}\label{skewsymmetry-condition-replicator}
\bar{Q_1}\,{\rm diag}(d_i)_i\begin{pmatrix}
a_{11}&a_{12}&\ldots&a_{1(n-1)}\\
a_{21}&a_{22}&\ldots&a_{2(n-1)}\\
\vdots&\vdots&\vdots\,\,\vdots\,\,\vdots&\vdots\\
a_{(n-1)1}&a_{(n-1)2}&\ldots&a_{(n-1)(n-1)}
\end{pmatrix}\,Q_2,
\end{equation}
is skew-symmetric, then the function
\[H_D(u)=\sum_{i=1}^{n-1}d_i(\sum_{k=1}^{n-1}\bar{q}_{ki}q_k) u_i+\sum_{i=1}^{n-1}d_i((\sum_{k=1}^{n-1}\bar{q}_{ki}q_k)-\bar{q}_{ii})e^{u_i}\]
is a constant of motion for the vector field $\tilde{X}_B(u)$ defined at~\eqref{jacobian-form2} and, consequently, the function 
\[H_d\circ\phi^{-1}(x)=H_{D}\circ\phi^{-1}(x)=\sum_{i=1}^{n-1}d_i(\sum_{k=1}^{n-1}\bar{q}_{ki}q_k) \log(\frac{x_i}{x_n})+\sum_{i=1}^{n-1}d_i((\sum_{k=1}^{n-1}\bar{q}_{ki}q_k)-\bar{q}_{ii})\frac{x_i}{x_n}\]
is a constant of motion for the replicator vector field $X_A$.
 \end{corollary}
\begin{remark}
The last line of the payoff matrix $A$ being null does not cause any loss of generality due to Item $1$ of Lemma~\ref{equivalence}. 
\end{remark}
\begin{proof} ({\bf Corollary~\ref{cor:replicator-general}})
The replicator equations with this type of payoff matrices have the point $q$ as a formal equilibrium point. Since $(\bar{Q_1})^tQ_1$ is diagonal the matrix
\[D=\bar{Q_1}{\rm diag}(d_i)_i,\]
satisfies Condition $(2)$ of Theorem~\ref{thm:consta-of-motion}. The matrix $B=(-1)EAE^t$ in Equation \eqref{jacobian-form2} is 
\begin{equation}\label{eq:from-B-to-A-general}
B=(-1)\begin{pmatrix}
a_{11}&a_{12}&\ldots&a_{1(n-1)}\\
a_{21}&a_{22}&\ldots&a_{2(n-1)}\\
\vdots&\vdots&\vdots\,\,\vdots\,\,\vdots&\vdots\\
a_{(n-1)1}&a_{(n-1)2}&\ldots&a_{(n-1)(n-1)}
\end{pmatrix}Q_2.\end{equation}
which yields
\[ DB=(-1)\bar{Q_1}\,{\rm diag}(d_i)_i\begin{pmatrix}
a_{11}&a_{12}&\ldots&a_{1(n-1)}\\
a_{21}&a_{22}&\ldots&a_{2(n-1)}\\
\vdots&\vdots&\vdots\,\,\vdots\,\,\vdots&\vdots\\
a_{(n-1)1}&a_{(n-1)2}&\ldots&a_{(n-1)(n-1)}
\end{pmatrix}\,Q_2.\]
The matrix at~\eqref{skewsymmetry-condition-replicator} being skew-symmetric means that $DB$ is skew-symmetric. Now, applying Theorem~\ref{thm:consta-of-motion} finishes the proof of the corollary. 
\end{proof}

\begin{remark}\label{D-for-internal}
Note that the matrix $D$ introduced in the proof of Corollary~\ref{cor:replicator-general} does not exhaust all possible matrices $D$ satisfying Condition $(2)$ of Theorem~\ref{thm:consta-of-motion}.  For example,  take the first column  $(d_{11},...,d_{1(n-1)})^t$ of $D$, then Condition $(2)$ implies that this column is a solution for the system
\[
\begin{pmatrix}
q_1&q_2&q_3&\ldots&q_{(n-1)}\\
q_1&q_2-1&q_3&\ldots&q_{(n-1)}\\
q_1&q_2&q_3-1&\ldots&q_{(n-1)}\\
\vdots&\vdots&\vdots&\vdots&\vdots\\
q_1&q_2&q_3&\ldots&q_{(n-1)}-1
\end{pmatrix}\begin{pmatrix}d_{11}\\d_{12}\\d_{13}\\\vdots\\d_{1(n-1)}\end{pmatrix}=\begin{pmatrix}0\\0\\0\\\vdots\\0\end{pmatrix},
\]

Simple calculation shows that $d_{1j}=d_{1(n-1)}$ for $j=2,\ldots,n-2$ and 
\[q_1d_{11}-(q_1+q_n)d_{1(n-1)}=0\]
where we used the fact that $q_2+q_3+\ldots+q_{n-1}-1=q_1+q_n$. This show that if $q_1=q_n=0$ then the first column of $D$ has two free variables. The same holds for any other column. This also justifies imposing the condition $q_n\neq0$ in Corollary~\ref{cor:replicator-general}.  
For an internal equilibrium point $q\in(\Delta^{n-1})^\circ$ we have
$$
D=\begin{pmatrix}\frac{q_1+q_n}{q_1}&1&\ldots&1\\
1&\frac{q_2+q_n}{q_2}&\ldots&1\\
\vdots&\vdots&&\vdots\\
1&1&\ldots&\frac{q_{n-1}+q_n}{q_{n-1}}
\end{pmatrix}{\rm diag}(d_i)_i,
$$
where $d_1,\ldots, d_{n-1}\in\mathbb{R}$ are to be determined by Condition $1$ of Theorem~\ref{thm:consta-of-motion}.
\end{remark}

\begin{example}\label{Ex:interior-equilibrium}
Let $q=(\frac{1}{n},\ldots,\frac{1}{n})$ then the payoff matrix of replicator equation which have $q$ as a formal equilibrium point are of the form
\[A=\begin{pmatrix}
a_{11}&a_{12}&\ldots&a_{1(n-1)}&-\sum_{j=1}^{n-1}a_{1j}\\
a_{21}&a_{22}&\ldots&a_{2(n-1)}&-\sum_{j=1}^{n-1}a_{2j}\\
\vdots&\vdots&\vdots\,\,\vdots\,\,\vdots&\vdots&\vdots\\
a_{(n-1)1}&a_{(n-1)2}&\ldots&a_{(n-1)(n-1)}&-\sum_{j=1}^{n-1}a_{(n-1)j}\\
0&0&\ldots&0&0
\end{pmatrix}.\]
By Remark~\ref{D-for-internal} 
\[\bar{Q_1}=\begin{pmatrix}2&1&\ldots&1\\
1&2&\ldots&1\\
\vdots&\vdots&&\vdots\\
1&1&\ldots&2
\end{pmatrix}\]
Clearly, $Q_1$ is invertible and  $Q_2=\frac{1}{n}\bar{Q}_1$.  Note that $\bar{Q}_1$ and $Q_2$ are symmetric then for any given matrix $G$
\begin{align*}
\bar{Q}_1 G Q_2&+ (Q_2)^t G^t (\bar{Q}_1)^t=0\,\Longleftrightarrow\, n(G+G^t)=(Q_2)^{-1}0(Q_2)^{-1}=0,
\end{align*}
which shows that the condition~\eqref{skewsymmetry-condition-replicator} is equivalent to 
\[{\rm diag}(d_i)_i\begin{pmatrix}
a_{11}&a_{12}&\ldots&a_{1(n-1)}\\
a_{21}&a_{22}&\ldots&a_{2(n-1)}\\
\vdots&\vdots&\vdots\,\,\vdots\,\,\vdots&\vdots\\
a_{(n-1)1}&a_{(n-1)2}&\ldots&a_{(n-1)(n-1)}
\end{pmatrix}\] being skew-symmetric. In this case the function
\[H_{D}\circ\phi^{-1}(x)=\sum_{i=1}^{n-1}d_i \left(\log(\frac{x_i}{x_n})-(\frac{x_i}{x_n})\right)\]
is the constant of motion for the replicator vector field $X_A$.
 \end{example}


\section{Hamiltonian Formalism} 
\label{dirac}
In this section, we first provide a simple introduction to (pre-)symplectic, poisson, Dirac and big-isotropic structures on $\mathbb{R}^m$ which is suitable for presenting our results on Hamiltonian formalism of replicator equations. Then we state the results on Hamiltonian formalism of replicator equations.

{\bf (pre-)Symplectic case:} Let $\omega$ be a closed two form on $\mathbb{R}^m$. It defines a linear vector bundle map 
$\omega^\sharp:T\mathbb{R}^m\to T^\ast\mathbb{R}^m$ by $X\mapsto \omega(X,.)$. We use the notation $\omega^\sharp(u)$ for the representing matrix, in the canonical basis, of the linear map $\omega^\sharp(u):T_u\mathbb{R}^m\to T_u^\ast\mathbb{R}^m$. If $\omega^\sharp(u)$ is invertible for every $u\in\mathbb{R}^m$ then $\omega$ is a \emph{symplectic structure} on $\mathbb{R}^m$ and $m$ is necessarily even. Relaxing the invertibility condition on $\omega^\sharp(u)$, the two form $\omega$ is called a \emph{presymplectic structure}. In both cases a Hamiltonian vector field $X_H$ is defined by $\omega^\sharp.X_H=d H$. 
In the sympletic case the flow of $X_H$ is volume preserving and as a consequence it has no attractor. Presymplectic structure does not put any restriction on the dimension of $\mathbb{R}^m$ but one loses the volume preserving property. Furthermore, in the presymplectic case, the vector field $X_H$ can be non-singular in the critical point of $H$. Since the level set of a critical point of $H$ is not a submanifold, there is more chance for chaotic behavior. In fact if $X_0$ is in the kernel of $\omega^\sharp$ then $X_H+X_0$ is still Hamiltonian with the same Hamiltonian function $H$. The vector field $X_0$ can behave unpredictably. 

{\bf Poisson case:} Let $\pi$ be a bivector on $\mathbb{R}^m$ i.e. a bilinear, antisymmetric map $\pi:T^\ast\mathbb{R}^m\times T^\ast\mathbb{R}^m\to\mathbb{R}$. Similar to the (pre)-symplectic case, it defines a linear vector bundle map $\pi^\sharp:T^\ast\mathbb{R}^m\to T\mathbb{R}^m$ by $\alpha\mapsto\pi(\alpha,.)$. We use the notation $\pi^\sharp(u)$ in the same manner as symplectic case. If $\pi$ satisfies the jacobi condition
\begin{equation}
\label{jacobi-condition-poisson} \sum_{l=1}^m (\pi_{ij}\frac{\partial \pi_{ik}}{\partial u_l}+\pi_{li}\frac{\partial \pi_{kj}}{\partial u_l}+\pi_{lk}\frac{\partial \pi_{jl}}{\partial u_l})=0\quad\forall i,j,k
\end{equation}
then it defines a \emph{poisson  structure} on $\mathbb{R}^m$ which is the same as defining a poisson  bracket i.e. a bilinear skew-symmetric bracket  $\{f,g\}:=(dg)^t\pi^\sharp df$ on $C^\infty(\mathbb{R}^m)$ which satisfies Leibniz's rule and Jacobi identity. A Hamiltonian vector field $X_H$ is defined by $X_H=\pi^\sharp dH$. 

  The Jacobi identity~\eqref{jacobi-condition-poisson} guaranties the integribilty of the distribution defined at every point $u\in\mathbb{R}^m$ by the image of the linear map $\pi^\sharp(u)$. 
Each leaf of this foliation have a symplectic structure induced by $\pi$. The dimension of the symplectic leaf passing through a given point is called the rank of the poisson structure at that point. The flow of $X_H$ preserves this foliation and its restriction to each one of these leafs is Hamiltonian in the symplectic sense. So in principle what one gets is a smooth bunch of Hamiltonian vector fields defined on the leaves of a symplectic foliation of $\mathbb{R}^m$. The Hamiltonian evolutionary games discussed in \cite{AD2014} are of this type. 

{\bf Dirac case:} Dirac structures were introduced in \cites{MR998124,MR951168} and have been used to study many systems such as systems with nonholonomic constraints \cite{MR2935374}, implicit Lagrangian and Hamiltonian systems \cite{MR2265464}, port-Hamiltonian systems \cite{MR2275732} and many more. Dirac structure unites and generalizes the poisson  and presymplectic structures (hence their "intersection" i.e. symplectic structure).

 We give the definition for general manifold $M$ before discussing the case $M=\mathbb{R}^m$ which is the case we need. Let $M$ be a manifold, then the vector bundle $\mathbb{T}M= TM\oplus T^\ast M$ is called the big tangent bundle or, in some literature, Pontryagin bundle. By $P_1:\mathbb{T}M\to TM$ and $P_2:\mathbb{T}M\to T^\ast M$ we, respectively, denote the
 projections on the first and second components. Denoting the natural pairing between vector field  $X\in\mathfrak{X}(M)$ and $1$-form $\alpha\in\Omega^1(M)$ by $\alpha(X)$, a natural pairing on the sections of $\mathbb{T}M$ is defined by
\begin{equation}
\label{pairing}
\ll (X,\alpha),(Y,\beta)\gg=\frac{1}{2}\left(\beta(X)+\alpha(Y)\right).
\end{equation} 
Let $L$ be a linear subbundle of $\mathbb{T}M$, its annihilator with respect to the pairing $\ll.,.\gg$ is defined as 
\[L^\perp:=\{(X,\alpha)\in \mathbb{T}M\,\,| \,\,\ll(X,\alpha),(Y,\beta)\gg=0\quad \forall (Y,\beta)\in L\}.\]
The pairing\, $\ll.,.\gg$ \, is neither positive definite nor negative definite. As a consequence for a given linear subbundle $L$ of   $\mathbb{T}M$ the intersection $L\cap L^\perp$ can be non-empty.  Having this in mind, a linear subbundle $L\subset \mathbb{T}M$ is called \emph{isotropic} if $L\subseteq L^\perp$. If $L=L^\perp$ then $L$ is called maximal isotropic. Maximal isotropy implies that the dimension of the fibers of $L$ is equal to the dimension of $M$.

\begin{defn}
A \emph{Dirac structure} on a manifold $M$ is a maximal isotropic linear subbunlde $L\subset TM\oplus T^\ast M$  such that for every $(X,\alpha),(Y,\beta)\in L$ 
 \begin{equation}\label{courant-bracket}
\left([X,Y], \mathcal{L}_X\beta-\mathcal{L}_Y\alpha+\frac{1}{2}d(\alpha(Y)-\beta(X))\right)\in L,
\end{equation}
where $[.,.]$ denotes the Lie bracket between vector fields and $\mathcal{L}$ stands for Lie derivative. 
The left hand side of~\eqref{courant-bracket} is called \emph{Courant bracket} of two sections $(X,\alpha),(Y,\beta)$.

 Furthermore, a vector field $X\in\mathfrak{X}(M)$  is called \emph{Hamiltonian} with respect to Dirac structure $L$ if there exist a Hamiltonian  $H\in C^\infty(M)$ such that $(X,d H)\in L.$
\end{defn}

 We also consider big-isotropic structure which is a generalization of Dirac structure. Up to our knowledge, not much has been done regarding the Hamiltonian systems with big-isotropic structures as underlying structure. In  \cite{MR2343378} the author studies the geometry of these structures and in \cite{MR2349409} he studies Hamiltonian systems in this context, providing some reduction theorems for this type of Hamiltonian systems. 

\begin{defn}  \emph{A big-isotropic structre} is an isotropic linear subbundle $L\subset TM\oplus T^\ast M$ which satisfies~\eqref{courant-bracket} and a vector field $X\in\mathfrak{X}(M)$  is called \emph{Hamiltonian} with respect to big-isotropic structure $L$ if there exist a Hamiltonian  $H\in C^\infty(M)$ such that $(X,d H)\in L.$
\end{defn}
   The following example shows that Dirac structure unifies symplectic and poisson  structures. 
\begin{example}\label{example:symplectic}
A (pre-)symplectic form $\omega$ and a poisson $\pi$ define the Dirac structures 
$L_\omega=\{(X,\omega^\sharp(X)) | X\in \mathfrak{X}(M)\}, $
respectively,
$L_\pi=\{(\pi^\sharp(\alpha),\alpha) | \alpha\in \Omega^1(M)\}.$

The skew-symmetricness of $\omega$ and $\pi$ yields the maximal isotropy condition and clossness of $\omega$, respectively, the Jacobi identity~\eqref{jacobi-condition-poisson} yield~\eqref{courant-bracket}.  
\end{example}

A consequence of the fact that  sections of a Dirac$\backslash$big-isotropic structure $L$ are closed with respect to Courant bracket~\eqref{courant-bracket} is the integrability of the (possibly singular) distribution $P_1(L)$. Every leaf $S$ of the foliation generated by $P_1(L)$ is equipped with the closed two form
\begin{equation}\label{presymplectic-on-leaves}
\omega_S(P_1(X,\alpha),P_1(Y,\beta))=\frac{1}{2}\left(\alpha(Y)-\beta(X)\right)\quad \forall P_1(X,\alpha),P_1(Y,\beta)\in TS,
\end{equation}
i.e. $P_1(L)$ integrates to a presymplectic foliation. 

Now, we proceed discussing the type of Dirac structures on $\mathbb{R}^m$ which we will be using. Let $B,D\in{\rm Mat}_{m\times m}(\mathbb{R})$ be two matrices and $f(u)\in C^\infty(\mathbb{R}^m)$ be a smooth function such that $f(u)\neq0\quad\forall u\in\mathbb{R}^m$. We define two matrix valued functions $\bf{B}:\mathbb{R}^m\to{\rm Mat}_{m\times m}(\mathbb{R})$ and $\bf{D}:\mathbb{R}^m\to{\rm Mat}_{m\times m}(\mathbb{R})$ by ${\bf{B}}(u)=f(u)B$ and ${\bf{D}}(u)=f(u)D$.

\begin{lemma}\label{lemma:our-dirac-structure}
Let $L_{(\bf{B},\bf{D}^t)}$ be the linear subbundle of  $T\mathbb{R}^m\oplus T^\ast \mathbb{R}^m$ which is defined by $L_{({\bf B}, {\bf D^t})}(u):=\{({\bf B}(u)z,{\bf D}^t(u)z)|\quad \forall z\in\mathbb{R}^m\}$ at every point $u\in\mathbb{R}^m$. Then subbundle $L_{({\bf B},{\bf D}^t)}$ is a big-isotropic structure if and only if
\begin{equation}\label{dirac-local-1}
(DB+B^tD^t)(u)=0
\end{equation}
and it is a Dirac structure if and only if it satisfies~\eqref{dirac-local-1} and 
\begin{equation}\label{dirac-local-2}
(\ker B)\cap(\ker  D^t)=0.
\end{equation}  
\end{lemma}
\begin{proof}
Condition~\eqref{dirac-local-1} is equivalent to being isotropic and Condition~\eqref{dirac-local-2} guaranties the maximality. So it remain to prove the integrability condition. The section of $L_{({\bf B},{\bf D}^t)}$ are of the form $(X,\alpha)$ such that
\[X(u)=f(u)B\eta(u)\quad\mbox{and}\quad \alpha=f(u)D^t\eta(u),\]
 where $\eta:\mathbb{R}^m\to\mathbb{R}^m$ is an smooth map. Let us set $f(u)=1$ for a moment. For given sections $(X,\alpha)=(B\eta_1,D^t\eta_1)$ and $(Y,\beta)=(B\eta_2,D^t\eta_2)$ we have
 \[[X,Y]=J_YX-J_XY=B(J_{\eta_2}B\eta_1-J_{\eta_1}B\eta_2),\]
 where $J_F$ denotes the Jacobian of a given map $F:\mathbb{R}^m\to\mathbb{R}^m$. Now, if we show that 
 \begin{equation}\label{eq:closness-of-courant-part-2}
\mathcal{L}_X\beta-\mathcal{L}_Y\alpha+\frac{1}{2}d(\alpha(Y)-\beta(X))=D^t(J_{\eta_2}B\eta_1-J_{\eta_1}B\eta_2),
 \end{equation}
 we are done for the case where $f(u)=1$. Equation~\eqref{eq:closness-of-courant-part-2} is a straightforward consequence of
\begin{itemize}
\item[-] The skew-symmetric property the matrix $DB$.
\item[-] Cartan magic formula $\mathcal{L}_X\beta=i_X d\beta+d(i_X \beta)$.
\item[-] The fact that $i_{B\eta_1} d(D^t\eta_2)=\left(D^tJ_{\eta_2}-(D^tJ_{\eta_2})^t\right)B\eta_1$ and similar equality for $i_{B\eta_2} d(D^t\eta_1)$
\item[-] The fact that $d<B\eta_1,D^t\eta_2>=(J_{\eta_1})^tB^tD^t\eta_2+(J_{\eta_2})^tBD\eta_1.$
 \end{itemize}
 For the general case with $f(u)$, with have 
 \begin{align*}
 &[X,Y]=J_YX-J_XY=f(u)B(J_{f(u)\eta_2}B\eta_1-J_{f(u)\eta_1}B\eta_2),\quad \mbox{and}\\
 &\mathcal{L}_X\beta-\mathcal{L}_Y\alpha+\frac{1}{2}d(\alpha(Y)-\beta(X))=f(u)D^t(J_{f(u)\eta_2}B\eta_1-J_{f(u)\eta_1}B\eta_2),
 \end{align*}
 which finishes the proof.
\end{proof}

Following Corollary is an immediate consequence of the fact that for a given non-singular  matrix $W$ we have
\[DB\quad\mbox{ is skew-symmetric}\quad\Leftrightarrow  \quad W^tDBW\quad \mbox{ is skew-symmetric}.\]

\begin{corollary}\label{invertible-cases} Considering the setting of Lemma~\ref{lemma:our-dirac-structure}, we have
\begin{itemize}
\item[i)] If\,\, $B$ is invertible then $L_{({\bf B,D^t})}=L_{(I, D^tB^{-1})}$ is the Dirac structure generated by constant presymplectic form $\omega^\sharp=D^tB^{-1}$.  
\item[ii)] If\,\, $D^t$ is invertible the $L_{({\bf B,D^t})}=L_{(B,(D^t)^{-1},I)}$ is the Dirac structure generated by constant poisson structure $\pi^\sharp= B(D^t)^{-1}$. 
\end{itemize}
\end{corollary}

Considering Dirac$\backslash$big-isotropic structure $L_{({\bf B},{\bf D}^t)}$, a pair $(X,d H)$ is a Hamiltonian system if and only if there exist a function $\eta:\mathbb{R}^m\to \mathbb{R}^m$ such that $X={\bf B}\eta$ and $dH={\bf D}^t\eta$. Now, we are ready to state the main result of this section. 
\begin{theorem}\label{thm:hamiltonian-formalism}
Given a replicator equation, $X_A$ with a formal equilibrium point $q$, let $\tilde{X}_B(u)=B\eta_q(u)$ be the vector field defined at \eqref{jacobian-form2} rewritten using notation setting~\eqref{eta-formula}.  Also, consider that there exist an $(n-1)\times(n-1)$ matrix $D$ which satisfies hypothesizes $1,2$ of Theorem~\ref{thm:consta-of-motion} and let 
\begin{equation}\label{H_D(u)}H_D(u)=\sum_{i=1}^{n-1}(\sum_{k=1}^{n-1}d_{ki}q_k) u_i+\sum_{i=1}^{n-1}((\sum_{k=1}^{n-1}d_{ki}q_k)-d_{ii})e^{u_i},
\end{equation}
 \begin{equation}\label{H_D(x)}H_{D}\circ\phi^{-1}(x)=\sum_{i=1}^{n-1}(\sum_{k=1}^{n-1}d_{ki}q_k) \log(\frac{x_i}{x_n})+\sum_{i=1}^{n-1}((\sum_{k=1}^{n-1}d_{ki}q_k)-d_{ii})\frac{x_i}{x_n},
 \end{equation}
as in Theorem~\ref{thm:consta-of-motion},  
\begin{equation}\label{B-bold-D-bold}
{\bf B}(u):=(1+\sum_{i=1}^{n-1}e^{u_i})B\quad\mbox{and}\quad {\bf D}^t(u):=(1+\sum_{i=1}^{n-1}e^{u_i})D^t \quad \forall u\in\mathbb{R}^m,
\end{equation}
and 
\begin{equation}\label{tilde-Y}
\tilde{Y}_B=({1+\sum_{i=1}^{n-1}e^{u_i}})\tilde{X}_B,
\end{equation}
then the pair $(\tilde{Y}_B, dH_D)$ is a Hamiltonian system with respect to Dirac $\backslash$ big-isotropic structure $L_{({\bf B,D^t})}$. 
\end{theorem}
\begin{proof} Note that by Lemma~\ref{lemma:our-dirac-structure} the subbundle $L_{({\bf B,D^t})}$ is a Dirac $\backslash$ big-isotropic structure. To prove the Theorem, we simply need to remind that $\tilde{X}_B=B\eta_q$ where $\eta_q$ is defined at~\eqref{eta-formula} and $$dH_D(u)=(1+\sum_{i=1}^{n-1}e^{u_i})D^t\eta_q,$$
i.e.
\[\tilde{Y}_B={\bf B}\eta\quad\mbox{and}\quad dH_D={\bf D^t}\eta,\]
so $(\tilde{Y}_B,d H_D)\in L_{({\bf B,D^t})}$. 
\end{proof}
\begin{corollary} Considering the setting of Theorem~\ref{thm:hamiltonian-formalism}, we have
\begin{itemize}
\item[1)] If the constant matrix $B$ is invertible then the vector field $\tilde{Y}_B$ is Hamiltonian with respect to constant (pre)symplectic structure $\omega^\sharp=D^tB^{-1}$ and with Hamiltonian $H_D(u)$ or equivalently, the vector field $\frac{1}{x_n}X_A$, restricted into the interior of simplex $\Delta^{n-1}$, is Hamiltonian with respect to presympletic structure $\omega^\sharp=(d (\phi)^{-1})^tD^tB^{-1}(d(\phi)^{-1})$, where $\phi$ is defined at~\eqref{phi-alpha}, and Hamiltonian $H_D\circ(\phi)^{-1}(x)$. 
 \item[2)] If the constant matrix $D^t$ is invertible then the vector field $\tilde{Y}_B$ is Hamiltonian with respect to constant poisson structure $\pi^\sharp=B(D^t)^{-1}$ and with Hamiltonian $H_D(u)$ or equivalently, the vector field $\frac{1}{x_n}X_A$, restricted into the interior of simplex $\Delta^{n-1}$, is Hamiltonian with respect to poisson structure $\pi^\sharp=(d\phi)B(D^t)^{-1}(d \phi)^t$ and Hamiltonian $H_D\circ(\phi)^{-1}(x)$. 
 \end{itemize} 
\end{corollary}

\section{Conservative Lotka-Volterra equations}
\label{lotka-volterra}

As mentioned in the introduction, the following  class of o.d.e.'s 
\begin{equation}\label{lv}
 \frac{d y_i}{dt} = y_i\left( \, r_i+\sum_{j=1}^m a^\prime_{ij}\,y_j\right) 
\; (1\leq i\leq m)\;,
\end{equation}
are known as Lotka-Volterra (LV) systems. The matrix $A^\prime=\matbrackets{a^\prime_{ij}}{ij}$ is called the system's interaction matrix. We will denote by $Y_{(A^\prime,r)}(y)$ the vector filed defined by Equations~\eqref{lv}. Following Theorem, \cite[Theorem 7.5.1]{HS}, shows that the replicator equation in $n$ variables $x_1,\ldots,x_n$ is equivalent to the Lotka-Volterra equation in $m=n-1$ variables $y_1,\ldots,y_{n-1}$. 

\begin{theorem}\label{thm:equivalence-lv-replicator}
The differentiable and invertible map from the interior of the simplex $\Delta^{n-1}$ onto $\mathbb{R}^{n-1}_+$ defined by 
\[y_i=\frac{x_i}{x_n}\quad i=1,\ldots,n-1\] 
 pushes the replicator vector field , $X_A(x)$, forward the time re-parametrization by $x_n$ of Lotka-Volterra vector field, $Y_{(A^\prime,r)}(y)$, where $r_i=a_{in}-a_{nn}$ and $a^\prime_{ij}=a_{ij}-a_{nj}.$ The inverse map is defined by
\[x_i=\frac{y_i}{1+\sum_{j=1}^{n-1}y_j},\quad i=1,\ldots,n-1\] 
and $x_n=\frac{1}{1+\sum_{j=1}^{n-1}y_j}$. 
\end{theorem}
\begin{proof}
See \cite[Theorem 7.5.1]{HS}.
\end{proof}

By Theorem~\ref{thm:equivalence-lv-replicator}, a LV system
$$ \frac{d y_i}{dt} = y_i\left( \, r_i+\sum_{j=1}^{n-1} a^\prime_{ij}\,y_j\right) 
\; (1\leq i\leq n-1)\;,$$ is equivalent to the replicator equation with payoff matrix 
\[A=\begin{pmatrix}
A^\prime&r\\0&0
\end{pmatrix}.
\]
 If $q^\prime\in\mathbb{R}^{n-1}_+$ is an equilibrium point for LV systems then
 $$q=(\frac{q^\prime_1}{1+\sum_{j=1}^{n-1}q^\prime_j},\ldots,\frac{q^\prime_{n-1}}{1+\sum_{j=1}^{n-1}q^\prime_j},\frac{1}{1+\sum_{j=1}^{n-1}q^\prime_j})$$ is an equilibrium for its equivalent replicator equation. In fact, $q$ is a formal equilibrium since $(Aq)_i=(Aq)_j=0$ for every $i,j$ and clearly 
 $\sum_{i=1}^n q_i=1$. Instead of considering $q^\prime\in\mathbb{R}^{n-1}_+$ we will consider $q^\prime\in\mathbb{R}^{n-1}$. 
 \begin{defn}\label{defn:formal-equilibrium-LV}
 A point $q^\prime\in\mathbb{R}^{n-1}$ is called a formal equilibrium for a given Lotka-Voltra equation if it is an equilibrium for the trivial extension of the given Lotka-Volterra equation to $\mathbb{R}^{n-1}$. 
 \end{defn}
 It is clear that a formal equilibrium of a given Lotka-Volterra equations yields a formal equilibrium for its equivalent replicator equation. Following Theorem is a consequence of Corollary~\ref{cor:replicator-general} and Theorem~\ref{thm:hamiltonian-formalism}.
 
 \begin{theorem}\label{thm:result-for-lv}
 Given a Lotka-Volterra equation, $Y_{(A^\prime,r)}(y)$, assume that it has a formal equilibrium $q^\prime$. Let $q=(\frac{q^\prime_1}{1+\sum_{j=1}^{n-1}q^\prime_j},\ldots,\frac{q^\prime_{n-1}}{1+\sum_{j=1}^{n-1}q^\prime_j},\frac{1}{1+\sum_{j=1}^{n-1}q^\prime_j})$ and $Q_1,Q_2$ the matrices defined in \eqref{Q1} and \eqref{Q2}.  If there exist an $(n-1)\times(n-1)$ matrix $D$ such that 
\begin{enumerate}
\item The matrix $DA^\prime Q_2$ is anti-symetric.  
\item The matrix $D^tQ_1$ is diagonal. 
\end{enumerate}
Then the function  
 $$H_{D}(y)=\sum_{i=1}^{n-1}(\sum_{k=1}^{n-1}d_{ki}q_k) \log y_i+\sum_{i=1}^{n-1}((\sum_{k=1}^{n-1}d_{ki}q_k)-d_{ii})y_i$$
  is a constant of motion  for  $Y_{(A^\prime,r)}(y)$. Furthermore, the Lotka-Volterra equation $Y_{(A^\prime,r)}(y)$ restricted to interior of $\mathbb{R}^{n-1}_+$ is equivalent to vector filed $\tilde{Y}_B(u)$ , defined at~\eqref{tilde-Y}, with $B=A^\prime Q_2$. This vector field is Hamiltonian having  $H_D(u)$, defined at~\eqref{H_D(u)}, as its Hamiltonian function and  Dirac $\backslash$ big-isotropic structure $L_{({\bf B,D^t})}$, where $\bf{B}$ and $\bf{D}$ are defined at \eqref{B-bold-D-bold}, as underlying geometric structure.
 \end{theorem}
 \begin{proof}
 Notice that the payoff matrix \[A=\begin{pmatrix}
A^\prime&r\\0&0
\end{pmatrix}
\]
is of the type which is considered in Corollary~\ref{cor:replicator-general}. As shown there for this type of payoff matrices the matrix $B=(-1)EAE^t$ in Equation~\eqref{jacobian-form2} is equal $A^\prime Q_2$. Considering this fact, the statement of Theorem is simply rewriting the statement of Corollary~\ref{cor:replicator-general} and Theorem~\ref{thm:hamiltonian-formalism}.
 \end{proof}
 
 Let us consider a Lotka-Volterra equation $Y_{(A^\prime,r)}(y)$  which is conservative in the sense that it has an equilibrium $q\in(\mathbb{R}^{(n-1)})^\circ$ and there exists diagonal matrix $D^\prime$  with strictly positive diagonal elements such that $A^\prime D^\prime$ is skew-symmetric.  Applying Item $2$ of Lemma~\ref{equivalence}, with $c_i=\frac{(1+\sum_{j=1}^{n-1}q^\prime_j)}{q^\prime_i}$ for $i=1,\ldots,n-1$ and $c_n=(1+\sum_{j=1}^{n-1}q^\prime_j)$, to the equivalent replicator equation $X_A(x)$ of $Y_{(A^\prime,r)}(x)$ shows that $Y_{(A^\prime,r)}(y)$ is  equivalent to time re-parametrization of replicator equation $X_{\tilde{A}}(\bar{x})$ where 
 \[\tilde{A}=\begin{pmatrix}
A^\prime&r\\0&0
\end{pmatrix}\, {\rm diag}(q^\prime_1,\ldots,q^\prime_{n-1},1),\]
and $\bar{x_i}=\frac{c_ix_i}{\sum_{j=1}^{n-1}c_jx_j}$ for $i=1,\ldots,n$.  The point $q=(\frac{1}{n},\ldots,\frac{1}{n})$ is an equilibrium point of the replicator equation $X_{\tilde{A}}(\bar{x})$. So by Example~\ref{Ex:interior-equilibrium} it is Hamiltonian in the setting of this paper if there exists $d\in\mathbb{R}^{n-1}$ such that 
 \begin{equation}\label{condition-of-example}
 {\rm diag}(d_i)_i A^\prime {\rm diag}(q^\prime_i)_i
\end{equation}
is skew-symmetric. Now, note that
\[A^\prime D^\prime+D^\prime (A^\prime)^t=0\,\Longrightarrow\,(D^\prime)^{-1}A^\prime+(A^\prime)^t(D^\prime)^{-1}=(D^\prime)^{-1}0(D^\prime)^{-1}=0,\]
i.e. $A^\prime D^\prime$ being skew-symmetric is equivalent to  
$(D^\prime)^{-1}A^\prime$ being skew-symmetric. Clearly, this implies that
\[{\rm diag}(q^\prime_1,\ldots,q^\prime_{n-1})
(D^\prime)^{-1}A^\prime{\rm diag}(q^\prime_1,\ldots,q^\prime_{n-1})\]
is skew-symmetirc. So  $d=(\frac{q^\prime_1}{ d^\prime_1},\ldots,\frac{q^\prime_{n-1}}{ d^\prime_{n-1}})_i$ satisfy condition~\eqref{condition-of-example}.  This shows that the conservative Lotka-Volterra equations introduced here include the ones obtained by gauge transformations. The Hamiltonian obtained by Example~\ref{Ex:interior-equilibrium} is
\[H_{D}\circ\phi^{-1}(\bar{x})=\sum_{i=1}^{n-1}\frac{q^\prime_i}{d^\prime_i} (\log(\frac{\bar{x}_i}{\bar{x}_n})-\frac{\bar{x}_i}{\bar{x}_n}).\]
a simple calculation together with Theorem~\ref{thm:equivalence-lv-replicator} shows that 
\[H_D(y)=\sum_{i=1}^{n-1}\frac{1}{ d^\prime_i}(q^\prime_i\log(y_i)-y_i)-\sum_{i=1}^{n-1}\frac{q^\prime_i}{ d^\prime_i}\log(q^\prime_i),\]
is the constant of motion obtained by our method for Lotka-Volterra equation $Y_{(A^\prime,r)}(y)$.  Note that 
\[H_D(y)=-H(y)-\sum_{i=1}^{n-1}\frac{q^\prime_i}{ d^\prime_i}\log(q^\prime_i),\]
where $H(y)$ is the constant of motion~\eqref{lv-constant-motion-intro}. So our method, basically, yields the same constant of motion. 

 On the other hand, let $Y_{(A^\prime,r)}(y)$ be a LV system  with interaction matrix
\[A^\prime=\begin{pmatrix}
0&\lambda_1&0&0\\
\lambda_2&0&0&0\\
0&0&\lambda_3&\lambda_5\\
0&0&\lambda_6&\lambda_4
\end{pmatrix}\]
where $\lambda_1\lambda_2\lambda_3\lambda_4\neq0$ and $r=-(\lambda_1,\lambda_2,\lambda_3+\lambda_5,\lambda_4+\lambda_6)$. The point $q^\prime=(1,1,1,1)$ is an equilibrium point of $Y_{(A^\prime,r)}(y)$. This implies that the point $q=(\frac{1}{5},\frac{1}{5},\frac{1}{5},\frac{1}{5},\frac{1}{5})$ is a formal equilibrium for the equivalent replicator equation. The matrix 
$${\rm diag}(\lambda_2,-\lambda_1,0,0)A^\prime,$$ is skew-symmetric. So $Y_{(A^\prime,r)}(y)$ is conservative in our setting. Note that $A^\prime D^\prime$ being skew-symmetric together with $D^\prime>0$ implies that $a^\prime_{ii}=0$ for $i=1,..,n-1$. Since $\lambda_3\lambda_4\neq 0$ it does not exist $D^\prime>0$ that makes $A^\prime D^\prime$ skew-symmetric. This shows that our approach enlarges the set of conservative LV systems.  

\begin{remark}\label{predator-prey}
In spite  of being a simple example one can see that considering matrix $A^\prime$ as interaction matrix between one predator specie $y_1$ and three prey species $y_2,y_3,y_4$ then the interaction between prey species $y_3,y_4$ is considered.  
\end{remark}

As for replicator equations, the equivalent replicator equations of $Y_{(A^\prime,r)}(y)$ has the payoff matrix
\[A=\begin{pmatrix}
0&\lambda_1&0&0&-\lambda_1\\
\lambda_2&0&0&0&\lambda_2\\
0&0&\lambda_3&\lambda_5&-\lambda_3-\lambda_5\\
0&0&\lambda_6&\lambda_4&-\lambda_4-\lambda_6\\
0&0&0&0&0
\end{pmatrix}.\]
This matrix can not be skew-symmetrized neither by  Item $1$ nor Item $2$ of Lemma~\ref{equivalence}.

\section*{Acknowledgments}

The author would like to thanks Pedro Duarte for the helpful discussions and comments. The preceding work \cite{AD2014} on this topic would not be possible without his collaboration. The author would also like to thanks Rui Loja Fernandes for helpful discussions and comments.


\end{document}